\documentclass{svproc}
\usepackage{url}

\usepackage{amssymb}
\usepackage{epsfig}
\usepackage{amsmath}
\usepackage[mathscr]{eucal}
\usepackage{subfigure}
\usepackage{graphicx}

\renewcommand{\qed}{\hfill{\ \ \rule{2mm}{2mm}} \vspace{0.2in}}

\begin{document}
\mainmatter              
\title{Size of local finite field Kakeya sets}
\titlerunning{Size of local finite field Kakeya sets}  
%
\author{Ghurumuruhan Ganesan\inst{1}}
%
\authorrunning{G. Ganesan} 
%
\tocauthor{G. Ganesan}
\institute{Institute of Mathematical Sciences, HBNI, Chennai\\
\email{gganesan82@gmail.com}}

\maketitle              

\begin{abstract}
Let~\(\mathbb{F}\) be a finite field consisting of~\(q\) elements and let~\(n \geq 1\) be an integer. In this paper, we study the size of local Kakeya sets with respect to subsets of~\(\mathbb{F}^{n}\) and obtain upper and lower bounds for the minimum size of a (local) Kakeya set with respect to an arbitrary set~\({\mathcal T} \subseteq \mathbb{F}^{n}.\)

\keywords{Local Kakeya Sets, Minimum Size, Probabilistic Method}
\end{abstract}
\renewcommand{\theequation}{\thesection.\arabic{equation}}
\setcounter{equation}{0}
\section{Introduction} \label{intro}

The study of finite field Kakeya sets is of interest from both theoretical and application perspectives. Letting~\(\mathbb{F}\) be a finite field containing~\(q\) elements and~\( n\geq 1\) be an integer, Wolff~\cite{wol} used counting arguments and planes to estimate that the minimum size of a global Kakeya set covering all vectors in~\(\mathbb{F}^{n}\) grows at least as~\(q^{n/2}.\) Later Dvir~\cite{deer} used polynomial methods to obtain sharper bounds (of the form~\(C\cdot q^{n}\)) on the minimum size of global Kakeya sets and for further improvements in the multiplicative constant~\(C,\) we refer to Saraf and Sudan~\cite{bhas}.

In this paper, we are interested in studying local Kakeya sets with respect to \emph{subsets} of~\(\mathbb{F}^{n}.\) Specifically, in Theorem~\ref{kak_thm}, we obtain upper and lower bounds for the minimum size of a Kakeya set with respect to a subset~\({\mathcal T} \subseteq \mathbb{F}^n.\)

The paper is organized as follows. In Section~\ref{pf_prop1}, we describe local Kakeya sets and state and prove our main result (Theorem~\ref{kak_thm}) regarding the minimum size of a local Kakeya set.

\setcounter{equation}{0}
\renewcommand\theequation{\thesection.\arabic{equation}}
\section{Local Kakeya sets}\label{pf_prop1}
Let~\(\mathbb{F}\) be a finite field containing~\(q\) elements and for~\(n \geq 1\) let~\(\mathbb{F}^{n}\) be the set of all~\(n-\)tuple vectors with entries belonging to~\(\mathbb{F}.\)

We say that a set~\({\mathcal K} \subseteq \mathbb{F}^{n}\) is a Kakeya set \emph{with respect to} the vector\\\(\mathbf{x}  = (x_1,\ldots,x_n) \in \mathbb{F}^{n}\) if there exists~\(\mathbf{y}  = \mathbf{y}(\mathbf{x}) \in \mathbb{F}^{n}\) such that the line
\begin{equation}\label{kak_def}
L(\mathbf{x},\mathbf{y}) := \bigcup_{a \in \mathbb{F}} \{\mathbf{y} + a \cdot \mathbf{x}\} \subseteq {\mathcal K},
\end{equation}
where~\(a \cdot \mathbf{x} := (ax_1,\ldots,ax_n).\) For a set~\({\mathcal T} \subseteq \mathbb{F}^{n},\) we say that~\({\mathcal K} \subseteq \mathbb{F}^{n}\) is a Kakeya set with respect to~\({\mathcal T}\) if~\({\mathcal K}\) is a Kakeya set with respect to every vector~\(\mathbf{x} \in {\mathcal T}.\)

The following result describes the minimum size of local Kakeya sets.
\begin{theorem}\label{kak_thm} Let~\({\mathcal T} \subseteq \mathbb{F}^{n}\) be any set with cardinality~\(\#{\mathcal T}\) an integer multiple of~\(q-1\) and let~\(\theta\left({\mathcal T}\right)\) be the minimum size of a Kakeya  set with respect to~\({\mathcal T}.\) We then have that
\begin{equation}\label{genh}
q\sqrt{M} + \min\left(0,q-\sqrt{M}\right) \leq \theta\left({\mathcal T}\right) \leq  q + q^{n} \left(1-\left(1-\frac{1}{q^{n-1}}\right)^{M-1}\right)
\end{equation}
where~\(M := \frac{\#{\mathcal T}}{q-1}.\)
\end{theorem}
For example suppose~\(M = \epsilon \cdot \left(\frac{q^{n}-1}{q-1}\right)\) for some~\(0 < \epsilon \leq 1.\) From the lower bound in~(\ref{genh}), we then get that~\(\theta\left({\mathcal T}\right)\) grows at least of the order of~\(q^{n/2}.\) Similarly, using the fact that~\(1-x \geq e^{-x-x^2}\) for~\(0 < x \leq \frac{1}{2},\) we get that
\[\left(1-\frac{1}{q^{n-1}}\right)^{M-1} \geq \exp\left(-\frac{M-1}{q^{n-1}} \left(1+\frac{1}{q^{n-1}}\right)\right) \geq e^{-\Delta},\]
where~\(\Delta := \frac{q\epsilon}{q-1} \left(1+\frac{1}{q^{n-1}}\right).\) From~(\ref{genh}) we then get that~\[\theta\left({\mathcal T}\right) \leq q + q^{n}(1-e^{-\Delta}).\]

In what follows we prove the lower bound and the upper bound in Theorem~\ref{kak_thm} in that order.

\subsection*{Proof of Lower Bound in Theorem~\ref{kak_thm}}
The proof of the lower bound consists of two steps. In the first step, we extract a subset~\({\mathcal N}\) of~\({\mathcal T}\) containing vectors that are non-equivalent. In the next step, we then use a high incidence counting argument similar to Wolff (1999) and estimate the number of vectors in a Kakeya set~\({\mathcal K}\) with respect to~\({\mathcal N}.\)

\emph{\underline{Step 1}}: Say that vectors~\(\mathbf{x}_1,\mathbf{x}_2 \in \mathbb{F}^{n}\) are \emph{equivalent} if~\(\mathbf{x}_1 = a \cdot \mathbf{x}_2\) for some~\(a \in \mathbb{F} \setminus \{0\}.\) We first extract a subset of vectors in~\({\mathcal T}\) that are pairwise non-equivalent. Pick a vector~\(\mathbf{x}_1 \in {\mathcal N}\) and throw away all the vectors in~\({\mathcal T},\) that are equivalent to~\(\mathbf{x}_1.\) Next, pick a vector~\(\mathbf{x}_2\) in the remaining set and again throw away the vectors that are equivalent to~\(\mathbf{x}_2.\) Since we throw away at most~\(q-1\) vectors in each step, after~\(r\) steps,
we are left with a set of size at least~\(\#{\mathcal T} - r(q-1).\) Thus the procedure continues for
\begin{equation}\label{m_def}
M = \frac{\#{\mathcal T}}{q-1}
\end{equation}
steps, assuming henceforth that~\(M\) is an integer.

Let~\({\mathcal N} =\{\mathbf{x}_1,\ldots,\mathbf{x}_M\} \subseteq {\mathcal T}\) be a set of size~\(M\) and let~\({\mathcal K}\) be a Kakeya set with respect to~\({\mathcal N},\) of minimum size. By definition (see~(\ref{kak_def})), there are vectors~\(\mathbf{y}_1,\ldots,\mathbf{y}_M\) in~\(\mathbb{F}^{n}\) such that the line~\(L(\mathbf{x}_i,\mathbf{y}_i) \subseteq {\mathcal K}\) for each~\(1 \leq i \leq M.\) Moreover, since~\({\mathcal K}\) is of minimum size,
we must have that
\begin{equation}\label{min_k}
{\mathcal K} = \bigcup_{j=1}^{M} \{L(\mathbf{x}_i,\mathbf{y}_i)\}.
\end{equation}

\emph{\underline{Step 2}}: To estimate the number of distinct vectors in~\({\mathcal K},\) suppose first that each vector in~\({\mathcal K}\) belongs to at most~\(t\) of the lines in~\(\{L(\mathbf{x}_i,\mathbf{y}_i)\}_{1 \leq i \leq M}.\) Since each line~\(L(\mathbf{x}_i,\mathbf{y}_i)\) contains~\(q\) vectors, we get that the total number of vectors in~\({\mathcal K}\) is bounded below by
\begin{equation}\label{first}
\#{\mathcal K} \geq \frac{Mq}{t}.
\end{equation}

Suppose now that there exists  a vector~\(\mathbf{v}\) in~\({\mathcal K}\) that belongs to at least~\(t+1\) lines~\(L_1,\ldots,L_{t+1} \subseteq \bigcup_{1 \leq i \leq M} \{L(\mathbf{x}_i,\mathbf{y}_i)\}.\) Because the vectors~\(\{\mathbf{x}_i\}_{1 \leq i \leq M}\) are not equivalent, any two lines~\(L(\mathbf{x}_i,\mathbf{y}_i)\) and~\(L(\mathbf{x}_j,\mathbf{y}_j)\) must have at most one point of intersection. To see this is true suppose there were scalars~\(a_1 \neq a_2\) and~\(b_1 \neq b_2\) in~\(\mathbb{F}\) such that
\[\mathbf{y}_1 + a_i \cdot \mathbf{x}_1 = \mathbf{y}_2 + b_i \cdot \mathbf{x}_2 \text{ for }i = 1,2.\]
Subtracting the equations we would then get~\((a_1-a_2) \cdot \mathbf{x}_1 = (b_1-b_2) \cdot \mathbf{x}_2,\)
contradicting the fact that~\(\mathbf{x}_1\) and~\(\mathbf{x}_2\) are not equivalent.

From the above paragraph, we get that any two lines in~\(\{L_i\}_{1 \leq i \leq t+1}\) have exactly one point of intersection, the vector~\(\mathbf{v}.\) Since each line~\(L_i\) contains~\(q\) vectors, the total number of vectors in~\(\{L_i\}_{1 \leq i \leq t+1}\) equals~\((q-1)(t+1)+1,\) all of which must be in~\({\mathcal K}.\) From~(\ref{first}) we therefore get that
\[\#{\mathcal K} \geq \min\left(\frac{Mq}{t} , (q-1)(t+1)+1\right)\]
and setting~\(t=\sqrt{M},\) we get
\[\#{\mathcal K} \geq q\sqrt{M} + \min\left(0,q-\sqrt{M}\right).\] From the expression for~\(M\) in~(\ref{m_def}), we then get~(\ref{genh}).~\(\qed\)

\subsection*{Proof of Upper Bound in Theorem~\ref{kak_thm}}
We use the probabilistic method. Let~\(\{\mathbf{x}_1,\ldots,\mathbf{x}_M\}\) be the set of non-equivalent vectors obtained in Step~\(1\) in the proof of the lower bound with\\\(M = \frac{\#{\mathcal T}}{q-1}\) (see~(\ref{m_def})). Let~\(\mathbf{Y}_1,\ldots,\mathbf{Y}_M\) be independently and uniformly randomly chosen from~\(\mathbb{F}^{n}\) and for~\(1 \leq i \leq M,\) set
\[{\mathcal S}_i := \bigcup_{j=1}^{i} \{L\left(\mathbf{x}_j,\mathbf{Y}_j\right)\},\]
where~\(L(\mathbf{x},\mathbf{y})\) is the line containing the vectors~\(\mathbf{x}\) and~\(\mathbf{y}\) as defined in~(\ref{kak_def}).

By construction, the set~\({\mathcal S}_M\) forms a Kakeya set with respect to~\({\mathcal T}.\) To estimate the expected size of~\({\mathcal S}_M,\) we use recursion. For~\(1 \leq i \leq M,\) let\\\(\theta_i := \mathbb{E}\#{\mathcal S}_i\) be the expected size of~\({\mathcal S}_i.\) Given~\({\mathcal S}_{i-1},\) the probability that a vector chosen from~\(\mathbb{F}^{n},\) uniformly randomly and independent of~\({\mathcal S}_{i-1},\) belongs to the set~\({\mathcal S}_{i-1}\) is given by~\(p_i := \frac{{\mathcal S}_{i-1}}{q^{n}}.\) Therefore
\[\mathbb{E}\#\left(L\left(\mathbf{x}_i,\mathbf{Y}_i\right) \bigcap {\mathcal S}_{i-1}\right) = q\cdot \mathbb{E}\left(\frac{\#{\mathcal S}_{i-1}}{q^{n}}\right)\]
and so
\begin{equation}\label{thet_rec}
\theta_i =\theta_{i-1} + q\left(1-\frac{\theta_{i-1}}{q^{n}}\right) = \theta_{i-1}\left(1-\frac{1}{q^{n-1}}\right) + q.
\end{equation}
Letting~\(a = 1-\frac{1}{q^{n-1}}\) and using~(\ref{thet_rec}) recursively, we get
\[\theta_i = a^{i-1}\cdot \theta_1+ q\cdot(1+a+\ldots+a^{i-2}) = a^{i-1} \cdot \theta_1 + \frac{q(1-a^{i-1})}{1-a}.\]
Using~\(\theta_1 = q,\) we then get that
\begin{eqnarray}
\theta_M  &=& a^{M-1} \cdot q + q^{n}\left(1-\left(1-\frac{1}{q^{n-1}}\right)^{M-1}\right) \nonumber\\
&\leq& q + q^{n} \left(1-\left(1-\frac{1}{q^{n-1}}\right)^{M-1}\right).\nonumber
\end{eqnarray}
This implies that there exists a Kakeya set with respect to~\({\mathcal T}\) of size at most~\(\theta_M.\)\(\qed\)

\subsection*{Acknowledgement}
I thank Professors V. Guruswami, V. Arvind, C. R. Subramanian and the referees for crucial comments that led to an improvement of the paper.
I also thank IMSc for my fellowships.


%
%

\end{document}